\documentclass[leqno,11pt]{article}

\usepackage{amssymb}
\usepackage{euscript}
\usepackage[dvips]{graphicx}
\usepackage{flafter}
\usepackage{pstricks}
\usepackage[latin1]{inputenc}
\usepackage{amsmath}
\usepackage{amsfonts}

\usepackage{mathrsfs} 

\evensidemargin\oddsidemargin

\begin{document}

\baselineskip=18pt \setcounter{page}{1}

\renewcommand{\theequation}{\thesection.\arabic{equation}}
\newtheorem{theorem}{Theorem}[section]
\newtheorem{acknowledgement}[theorem]{Acknowledgement}
\newtheorem{algorithm}[theorem]{Algorithm}
\newtheorem{axiom}[theorem]{Axiom}
\newtheorem{case}[theorem]{Case}
\newtheorem{claim}[theorem]{Claim}
\newtheorem{conclusion}[theorem]{Conclusion}
\newtheorem{condition}[theorem]{Condition}
\newtheorem{conjecture}[theorem]{Conjecture}
\newtheorem{corollary}[theorem]{Corollary}
\newtheorem{criterion}[theorem]{Criterion}
\newtheorem{definition}[theorem]{Definition}
\newtheorem{example}[theorem]{Example}
\newtheorem{exercise}[theorem]{Exercise}
\newtheorem{lemma}[theorem]{Lemma}
\newtheorem{notation}[theorem]{Notation}
\newtheorem{problem}[theorem]{Problem}
\newtheorem{proposition}[theorem]{Proposition}
\newtheorem{remark}[theorem]{Remark}
\newtheorem{solution}[theorem]{Solution}
\newtheorem{summary}[theorem]{Summary}
\newenvironment{proof}[1][Proof]{\noindent\textbf{#1.} }{\ \rule{0.5em}{0.5em}}
\newcommand{\eqnsection}{
\renewcommand{\theequation}{\thesection.\arabic{equation}}
    \makeatletter
    \csname  @addtoreset\endcsname{equation}{section}
    \makeatother}
\eqnsection

\def\r{{\mathbb R}}
\def\e{{\mathbb E}}
\def\p{{\mathbb P}}
\def\z{{\mathbb Z}}


\vglue20pt

\centerline{\Large\bf Branching random walk in  $\z^{4}$ with
branching}

\centerline{\Large\bf at the origin only}

\bigskip

\centerline{by}

\medskip

\centerline{Yueyun Hu,  Vladimir Vatutin\footnote{Supported in part
by RFBR grant 08-01-00078,
 and the programm ``Mathematical Control Theory'' RAS},   and Valentin  Topchii\footnote{Supported in part by the programm
 ``Contemporary Problems of Theoretical Mathematics of the MI SB RAS''}}

\medskip

\centerline{\it Universit\'e Paris XIII,    Steklov Mathematical
Institute and Omsk Branch of Sovolev Institute of Mathematics}

\bigskip
\centerline{\today}
\bigskip
\bigskip

{\leftskip=2truecm \rightskip=2truecm \baselineskip=15pt \small

\noindent{\slshape\bfseries Summary.} For the critical branching
random walk in  $\mathbb{Z}^{4}$ with branching at the origin only we find
the asymptotic behavior of the probability of the event that there
are particles at the origin at moment $t\rightarrow \infty $ and
prove a Yaglom type conditional limit theorem for the number of
individuals at the origin given that there are particles at the
origin.

\bigskip




} 

\bigskip
\bigskip

\section{Introduction}

Consider the following modification of a standard branching random
walk on $\mathbb{Z}^{d}$. The population is initiated at time $t=0$ by a
single particle. Being outside the origin the particle performs a continuous
time random walk on $\mathbb{Z}^{d}$ with infinitesimal transition matrix
\begin{equation*}
A=|a(\mathbf{x},\mathbf{y})|_{\mathbf{x},\mathbf{y}\in \mathbb{Z}%
^{d}},\qquad a(\mathbf{0},\mathbf{0})<0,
\end{equation*}%
until the moment when it hits the origin (that is, the time which
the particle spends at a point $\mathbf{x}\neq \mathbf{0}$ is
exponentially distributed with parameter
$a:=-a(\mathbf{0},\mathbf{0})$). At the origin it spends an
exponentially distributed time with parameter $1$ and then either
jumps to a point $\mathbf{y}\neq \mathbf{0}$ with probability
\begin{equation}
-(1-\alpha )a(\mathbf{0},\mathbf{y})a^{-1}(\mathbf{0},\mathbf{0})=:(1-\alpha
)\pi _{\mathbf{y}},  \label{Defpy}
\end{equation}%
or dies with probability $\alpha $ producing just before the death a random
number of children $\xi $ in accordance with offspring generating function

\begin{equation*}
f(s):=\mathbf{E}s^{\,\xi }=\sum_{k=0}^{\infty }f_{k}s^{k}, \qquad
0\le s\le 1.
\end{equation*}

At the birth moment the newborn particles are located at the origin and from
this point they begin their own branching random walks behaving
independently and stochastically the same as the parent individual.

This model was investigated in  \cite {AB1, AB2, BY1, BY2, Y2008}
where some basic equations for the probability generating
functions of the number of particles $\eta(t;\mathbf{x})$ at point
$\mathbf{x}\in \mathbb{Z}^{d}$ at time $t$ were deduced and, under
certain conditions, the asymptotic behavior of the moments of
$\eta(t;\mathbf{x})$ as $t\to\infty$ was investigated. For the
continuous counterpart of the above model, we refer to \cite{F94},
\cite{DF}, \cite{FLG95}, \cite{D96}  and the references therein.

In the present paper we assume that the branching random walk is initiated at time $t=0$ by a single individual
located at the origin and impose the following restrictions on the
characteristics of the process:

\textbf{Hypothesis (I)}: The underlying random walk on $\mathbb{Z}^{d}$ is
symmetric, irreducible and homogeneous
\begin{equation*}
a(\mathbf{x},\mathbf{y})=a(\mathbf{y},\mathbf{x}), \qquad a(\mathbf{x},\mathbf{y})=a(%
\mathbf{0},\mathbf{y}-\mathbf{x})=:a(\mathbf{y}-\mathbf{x}),
\end{equation*}%
where $a(\mathbf{x})\geq 0$ for $\mathbf{x}\neq \mathbf{0}$ and $a(\mathbf{0}%
)<0$; besides we assume that
\begin{equation}
\sum_{\mathbf{x}\in \mathbb{Z}^{d}}a(\mathbf{x})=0,\qquad b^{2}:=\sum_{%
\mathbf{x}\in \mathbb{Z}^{d}}|\mathbf{x}|^{2}a(\mathbf{x})<\infty ,
\label{var2}
\end{equation}%
where $|\mathbf{x}|^{2}:=\sum_{i=0}^{d}x_i^{2}$ for $\mathbf{x}%
=(x_1,...,x_i,...,x_d)\in \mathbb{Z}^{d}$.

Denote $h_{d}$ the probability that a particle, leaving the origin and
performing a random walk on $\mathbb{Z}^{d}$ satisfying Hypothesis (I), will
never come back. Observe that $h_{1}=h_{2}=0$, while $h_{d}\in (0,1),\,d\geq
3.$

\textbf{Hypothesis (II)}: The offspring process is critical
\begin{equation*}
\alpha f^{\prime}(1)+(1-\alpha )(1-h_{d})=1
\end{equation*}%
and $f^{(2)}(1)\in (0,\infty )$.

Here and in what follows for $k=2,3,\ldots$ we use the notation  $f^{(k)}(s):=d^kf(s)/ds^k.$

Clearly, for the dimensions $d=1,2$ the introduced criticality of the
branching random walk is reduced to the criticality of the offspring
generating function at the origin, that is to the condition $f^{\prime}(1)=1$. For the dimensions $d\geq 3$ this is not
the case.

Let $\mu _{0}(t)$ denote the number of particles in the process located at
time $t$ at the origin, $\mu (t)$ denote the number of particles in the
process at time $t$ outside the origin, and let $\eta (t)=\mu _{0}(t)+\mu
(t) $ be the total number of individuals at the process at moment $t$.

For the case $d=1$ the authors of papers \cite{TVY, TV03, TV04, TV04tv}
investigate the asymptotic properties of the probabilities $\left\{ \mu
_{0}(t)>0\right\} $ and $\left\{ \mu (t)>0\right\} $ and prove conditional
limit theorems for a properly scaled vectors $(\mu _{0}(t),\mu (t))$ by
studying an auxiliary branching process.
A modification of this model when  the initial number of particles at the origin is large
 was considered for the case $d=1$ in \cite{VX}.

The asymptotic behavior of the
probability $q(t):=\mathbf{P}\left( \mu _{0}(t)>0\right) $ for the dimensions
\thinspace $d\neq 4$ (again, by constructing an auxiliary branching process)
has been found in \cite{VT10}. However, the behavior of the probability $%
q(t) $ for the case $d=4$ remained an open problem. The present paper fills
this gap.

Our main result looks as follows.
\begin{theorem}
\label{Tonedim4Gen}For a branching random walk evolving in
$\mathbb{Z}^{4}$ and satisfying Hypotheses (I) and (II),
\begin{equation}
q(t)\sim C\frac{\log t}{t},  \qquad t\rightarrow \infty,
\label{AAsympq}
\end{equation}%
with  $C= 3 a (1-\alpha) \gamma_4 h_4^2 /( \alpha f^{(2)}(1))>0$ and
$\gamma_4$ a constant given by Lemma \ref{LTransit} below, and
\begin{equation}
\lim_{t\rightarrow \infty }\mathbf{P}\left( \frac{\mu _{0}(t)}{\mathbf{E}%
\left[ \mu _{0}(t)|\mu _{0}(t)>0\right] }\leq x\left\vert \mu
_{0}(t)>0\right. \right) =\frac{1}{3}+\frac{2}{3}\left( 1-e^{-2x/3}\right) ,\
x>0.  \label{CondLimit}
\end{equation}
\end{theorem}

\begin{remark}
It will be shown later on that
\begin{equation*}
\mathbf{E}\left[ \mu _{0}(t)|\mu _{0}(t)>0\right] =\frac{\mathbf{E}\mu
_{0}(t)}{\mathbf{P}\left( \mu _{0}(t)>0\right) }\sim  {3 \over \alpha f^{(2)}(1)  C^2 } \, \frac{t}{\log ^{2}t%
},\qquad  t\rightarrow \infty .
\end{equation*}
\end{remark}

To obtain Theorem \ref{Tonedim4Gen}, the crucial step is to show the
asymptotic for  the survival probability  $q(t)$, which satisfies
some convolution equation (see (\ref{Eqvq}) below).  It turns out
that a first order analysis of this equation only gives a rough
upper bound for $q(t)$ (Lemma \ref{Lrough}), and we need a second
order argument to get the exact asymptotic (Proposition
\ref{Lprecise}).

The rest of this paper is organized as follows: In Section
\ref{s:pre}, we recall some known facts and present some basic
evolution equations. The proof of Theorem \ref{Tonedim4Gen} will be
given in Section \ref{s:proofmain}.

\section{Auxiliary results and basic equations}\label{s:pre}

For further references we first list some results obtained in
\cite{VT10}. We forget for a wile that we deal with a branching random walk in $\mathbb{Z}%
^{d},d\geq 3,$ and consider only the motion of a particle performing
a random walk satisfying Hypothesis (I) without branching.

Denote
\begin{equation*}
p(t;\mathbf{x},\mathbf{y})=p(t;\mathbf{0},\mathbf{y}-\mathbf{x})=:p(t;%
\mathbf{y}-\mathbf{x})
\end{equation*}%
the probability that a particle located at moment $t=0$ at point $\mathbf{x}$
will be at point $\mathbf{y}$ at moment $t$.

\begin{lemma}
\label{LTransit}(see \cite{VT10}) Under Hypothesis \textrm{(I)} as $%
t\rightarrow \infty $%
\begin{equation*}
p(t;\mathbf{0})\sim \gamma _{d}\, t^{-d/2},\quad \gamma
_{d}>0\text{.}
\end{equation*}
\end{lemma}

Consider a particle located at the origin at moment $t=0$ and let
$\tau _{0} $ be the
time the particle spends at the origin. Denote by $\tau _{2,d}$ the
time needed for a particle which has left the origin to come back to the origin. Let  $\theta (\mathbf{y})$ be the time
needed for a particle located at point $\mathbf{y}\neq \mathbf{0}$
at time $0$ to hit the origin for the first time. Let
\begin{equation*}
G_{1}(t):=\mathbf{P}(\tau _{0}\leq t)=1-e^{-t},\ G_{2,d}(t):=\mathbf{P}
(\tau _{2,d}\leq t),\ G^{\mathbf{y}}(t):=\mathbf{P}\left( \theta (\mathbf{y}%
)\leq t\right)
\end{equation*}%
and $G_{2,d}:=\mathbf{P}(\tau _{2,d}<\infty )$. It follows from (\ref{Defpy}%
) that%
\begin{equation*}
G_{2,d}(t)=\sum_{\mathbf{y}\neq \mathbf{0}}\pi _{\mathbf{y}}G^{\mathbf{y}%
}(t).
\end{equation*}

\begin{lemma}
\label{l6} (see \cite{VT10}) Under Hypothesis \textrm{(I)}, $\tau
_{2,d},d\geq 3,$ is an improper random variable with
\begin{equation*}
\mathbf{P}(\tau _{2,d}=\infty )=h_{d}=1-G_{2,d}(\infty )=1-G_{2,d}.
\end{equation*}%
Besides, for $G_{2}(t):=G_{2,d}(t)/(1-h_{d})$ as $t\rightarrow \infty $
\begin{equation*}
1-G_{2}(t)\sim \frac{{2a\gamma _{d}h_{d}^{2}}}{\left( 1-h_{d}\right) {(d-2)}}%
t^{-d/2+1}.
\end{equation*}
\end{lemma}

Let
\begin{equation*}
K_{d}(t)=K(t)=:\alpha f^{\prime}(1)G_{1}(t)+(1-\alpha )(1-h_{d})\,G_{1}\ast
G_{2}(t).
\end{equation*}

\begin{lemma}
\label{Cden}(see \cite{VT10}) Under Hypothesis \textrm{(I) }for\textrm{\ }$%
d\geq 3$ as $t\rightarrow \infty $
\begin{eqnarray}
1-K_{d}(t) &\sim &(1-\alpha )\left( G_{2,d}-G_{2,d}(t)\right) \sim \frac{%
2a(1-\alpha )\gamma _{d}h_{d}^{2}}{d-2}\,t^{-d/2+1},  \label{Assstail} \\
k_{d}(t) &:=&K_{d}^{\prime }(t)\sim a(1-\alpha )\gamma
_{d}h_{d}^{2}\,t^{-d/2}.  \label{AsDD2}
\end{eqnarray}
\end{lemma}


Let
\begin{equation*}
F_{\mathbf{y}}(t;s_{1},s_{2}):=\mathbf{E}_{\mathbf{y}}s_{1}^{\mu
_{0}(t)}s_{2}^{\mu _{+}(t)},  \qquad 0\le s_1, s_2 \le1 ,
\end{equation*}%
be probability generating function for the number of particles at the origin and outside the origin at moment $t
\ge0$ in the branching random walk   generated by a single particle located at
point $\mathbf{y}$ at
moment $0$. By the total probability formula we have%
\begin{eqnarray*}
F_{\mathbf{0}}(t;s_{1},s_{2}) &=&s_{1}(1-G_{1}(t))+\int_{0}^{t}\alpha f(F_{%
\mathbf{0}}(t-u;s_{1},s_{2}))\,dG_{1}(u) \\
&&+(1-\alpha )\int_{0}^{t}\left( \sum_{\mathbf{y}\neq \mathbf{0}}\pi _{%
\mathbf{y}}F_{\mathbf{y}}(t-u;s_{1},s_{2})\right) dG_{1}(u),
\end{eqnarray*}%
while for $\mathbf{y}\neq \mathbf{0}$%
\begin{equation*}
F_{\mathbf{y}}(t;s_{1},s_{2})=s_{2}(1-G^{\mathbf{y}}(t))+\int_{0}^{t}F_{%
\mathbf{0}}(t-u;s_{1},s_{2})\,dG^{\mathbf{y}}(u).
\end{equation*}%
Hence%
\begin{eqnarray*}
F_{\mathbf{0}}(t;s_{1},s_{2}) &=& s_{1}(1-G_{1}(t))+\int_{0}^{t}\alpha f(F_{%
\mathbf{0}}(t-u;s_{1},s_{2}))\,dG_{1}(u) \\
&&+(1-\alpha )s_{2}\int_{0}^{t}\left( \sum_{\mathbf{y}\neq \mathbf{0}}\pi _{%
\mathbf{y}}(1-G^{\mathbf{y}}(t-u))\right) dG_{1}(u) \\
&&+(1-\alpha )\int_{0}^{t}\left( \sum_{\mathbf{y}\neq \mathbf{0}}\pi _{%
\mathbf{y}}\int_{0}^{t-v}F_{\mathbf{0}}(t-u-v;s_{1},s_{2})\,dG^{\mathbf{y}%
}(v)\right) dG_{1}(u) \\
&=&s_{1}(1-G_{1}(t))+\int_{0}^{t}\alpha f(F_{\mathbf{0}}(t-u;s_{1},s_{2}))%
\,dG_{1}(u) \\
&&+(1-\alpha )s_{2}\int_{0}^{t}{}(1-G_{2,d}(t-u))dG_{1}(u) \\
&&+(1-\alpha )\int_{0}^{t}F_{\mathbf{0}}(t-u;s_{1},s_{2})d\left( G_{1}\ast
G_{2,d}(u)\right) .
\end{eqnarray*}%
Using this relation and setting $F(t;s):=F_{\mathbf{0}}(t;s,1)=\mathbf{E}_{%
\mathbf{0}} s^{\mu _{0}(t)}$, we get that for all $0\le s\le 1$,
\begin{eqnarray}
F(t;s) &=&s(1-G_{1}(t))+(1-\alpha )(1-h_{d})(1-G_{2}(\cdot ))\ast G_{1}(t)
\notag \\
&+&\int_{0}^{t}\alpha f(F(t-u;s))\,dG_{1}(u)+(1-\alpha )h_{d}G_{1}(t)  \notag
\\
&+&\int_{0}^{t}(1-\alpha )(1-h_{d})F(t-u;s)\,d(G_{1}\ast G_{2}(u)).
\label{ur20}
\end{eqnarray}

Hence, letting $q(t;s):=1-F(t;s),$
\begin{equation}
\Phi (x):=f^{\prime}(1)x-(1-f(1-x))=:x\Psi (x)\sim \frac{f^{(2)}(1)}{2}%
x^{2}, \qquad x\rightarrow 0,  \label{AsymPhi}
\end{equation}
we deduce that for all $0\le s\le 1$,
\begin{equation}
q(t;s)=(1-s)(1-G_{1}(t))+q(\cdot ;s)\ast K(t)-\alpha \Phi (q(\cdot
;s))\ast G_{1}(t).  \label{EqGener}
\end{equation}%

Define  $q(t):=\mathbf{P}\left( \mu _{0}(t)>0\right) =q(t;0)$. We
have
\begin{equation*}
q(t)=1-G_{1}(t)+q\ast K(t)-\alpha \Phi (q(\cdot ))\ast G_{1}(t),
\qquad t\ge0.
\end{equation*}%
Note that%
\begin{equation*}
k_{d}(t)=\alpha f^{\prime}(1)e^{-t}+(1-\alpha )(1-h_{d})(G_{1}\ast
G_{2}(t))^{\prime }.
\end{equation*}%
Thus, we have%
\begin{equation}
q(t)=1-G_{1}(t)+\int_{0}^{t}q(t-u)\left( k_{d}(u)-\alpha \Psi
(q(t-u))e^{-u}\right) du.  \label{Eqvq}
\end{equation}

The previous arguments hold for any dimension $d$. Now we concentrate on the
case $d=4$ and recall that by (\ref{AsDD2})%
\begin{equation}
k_{4}(t)\sim c_{4}t^{-2}  \label{DenK}
\end{equation}%
as $t\rightarrow \infty $ where $c_{4}:=a\left( 1-\alpha \right) \gamma
_{4}h_{4}^{2}>0$. This asymptotic formula allows us to prove the following
statement.

\begin{lemma}  For any fixed $\varepsilon \in (0,1)$ and $p\in
(0,1]$, we have for $t \to \infty$
\begin{eqnarray}
I:&=&\int_{0}^{t\varepsilon }\frac{\log ^{p}\left( t-u+1\right)
}{t-u+1}  k_{4}(u)du=\frac{\log ^{p}(t+1)}{t+1}+\frac{c_{4}\log
^{1+p}t}{t^{2}} (1+o(1)), \label{Lintegr}  \\
I^{\ast }&:=&\int_{t\varepsilon }^{t}\frac{\log ^{p}\left( t-u+1\right) }{t-u+1%
}k_{4}(u)du=\frac{c_{4}\log ^{1+p}t}{\left( 1+p\right) t^{2}}+o\left( \frac{%
\log ^{1+p}t}{t^{2}}\right), \label{LintStar} \\
&& \int_{0}^{t}\frac{\log ^{p}\left( t-u+1\right)
}{t-u+1}k_{4}(u)du=\frac{\log ^{p}(t+1)}{t+1}+c_{4}
\frac{2+p}{1+p}(1+o(1)) \frac{\log ^{1+p}t}{t^{2}}. \label{CorInt}
\end{eqnarray}
\end{lemma}

\noindent\textbf{Proof}.  Clearly (\ref{CorInt}) follows from
(\ref{Lintegr}) and (\ref{LintStar}).  To  prove (\ref{Lintegr}),
we use the expansions (valid for $0\leq u\leq \varepsilon t$)
\begin{eqnarray*}
\log ^{p}\left( t-u+1\right) &=&\left( \log \left( t+1\right) -\frac{u}{t}%
+O\left( \frac{u^{2}+u}{t^{2}}\right) \right) ^{p} \\
&=&\log ^{p}(t+1)\left( 1-\frac{pu}{t\log \left( t+1\right) }+O\left( \frac{%
u^{2}+u}{t^{2}\log \left( t+1\right) }\right) \right)
\end{eqnarray*}%
and%
\begin{equation*}
\frac{1}{t-u+1}=\frac{1}{t+1}\left( 1+\frac{u}{t}+O\left( \frac{u^{2}+1}{t^{2}}%
\right) \right)
\end{equation*}%
implying%
\begin{eqnarray*}
&&\frac{\log ^{p}\left( t-u+1\right) }{t-u+1} \\
&=&\frac{\log ^{p}(t+1)}{t+1}+\frac{u\log ^{p}t}{t^{2}} +O\left( \frac{%
(u^{2}+u)\log ^{p}t}{t^{3}}\right) -\frac{pu}{t^{2}\log ^{1-p}(t+1)}+O\left(
\frac{u^{2}+u}{t^{3}\log ^{1-p}(t+1)}\right) .
\end{eqnarray*}%
Hence we have%
\begin{eqnarray*}
I &=&\frac{\log ^{p}(t+1)}{t+1}\int_{0}^{t\varepsilon }k_{4}(u)du+\frac{\log
^{p}t}{t^{2}}\int_{0}^{t\varepsilon }uk_{4}(u)du \\
&&+O\left( \frac{\log ^{p}t}{t^{3}}\right) \int_{0}^{t\varepsilon
}(u^{2}+1)k_{4}(u)du+O\left( \frac{\log ^{p-1}(t+1)}{t^{2}}\right)
\int_{0}^{t\varepsilon }uk_{4}(u)du \\
&=&\frac{\log ^{p}(t+1)}{t+1}+\frac{  c_{4}\log
^{1+p}t}{t^{2}}(1+o(1))+O\left( \frac{\log ^{p}t}{t^{2}}\right)
 \\
&=&\frac{\log ^{p}(t+1)}{t+1}+\frac{c_{4}\log ^{1+p}t}{t^{2}}(1+o(1)),
\end{eqnarray*}%
proving (\ref{Lintegr}). It remains to show (\ref{LintStar}).
Observe that
\begin{equation*}
I^{\ast }:=c_{4}(1+o(1))\int_{t\varepsilon }^{t}\frac{\log ^{p}\left(
t-u+1\right) }{t-u+1}\frac{1}{u^{2}}du.
\end{equation*}%
Further,%
\begin{eqnarray*}
\int_{t\varepsilon }^{t}\frac{\log ^{p}\left( t-u+1\right) }{t-u+1}\frac{1}{%
u^{2}}du &=&\int_{0}^{t(1-\varepsilon )}\frac{\log ^{p}\left( u+1\right) }{%
u+1}\frac{1}{\left( t-u\right) ^{2}}du \\
&=&\frac{1}{t^{2}}\int_{0}^{t(1-\varepsilon )}\frac{\log ^{p}\left(
u+1\right) }{u+1}du+O\left( \frac{1}{t^{3}}\int_{0}^{t(1-\varepsilon )}\log
^{p}\left( u+1\right) du\right) \\
&=&\left(\frac{1}{1+p}+ o(1)\right)\, \frac{\log ^{1+p}t}{t^{2}},
\end{eqnarray*} proving (\ref{LintStar}) and the Lemma. $\Box$

Let
\begin{equation*}
I_{k}(t):=\int_{2}^{t-2}\frac{\left( t-u\right) ^{k-1}\delta_1(t-u)\delta_2(u)}{\log ^{2k}(t-u)}%
\frac{du}{\log u},\qquad k=1,2,...,
\end{equation*}%
where $\delta_1(t)$ and $\delta_2(t)$ are bounded functions such that $\delta_1(t)\sim 1$ and
$\delta_2(t)\sim 1$ as $t\to\infty$.

\noindent We end this section by the following estimate on $I_k(t)$:

\begin{lemma}
\label{LbechInteg} For any fixed integer $k\ge1$, we have as $t\rightarrow \infty $%
\begin{equation*}
I_{k}(t)= \frac{1+o(1)}{k}\, \frac{t^{k}}{\log ^{2k+1}t}.
\end{equation*}
\end{lemma}

\noindent\textbf{Proof}. For any $\varepsilon \in (0, 1), $ if
$0\leq u\leq
\varepsilon t$ then%
\begin{equation*}
\log (t-u)=\log t+O(1)
\end{equation*}%
while if $t\geq u\geq \varepsilon t$ then%
\begin{equation*}
\log u=\log t+O(1).
\end{equation*}%
Clearly,%
\begin{equation*}
I_{1}(t)=\frac{\left( 1+o(1)\right) }{\log ^{2}t}\int_{2}^{t/2}\frac{1}{\log
u}du+\frac{\left( 1+o(1)\right) }{\log t}\int_{0}^{t/2}\frac{1}{\log ^{2}u}%
du= (1+o(1)) \frac{t}{\log ^{3}t}.
\end{equation*}%
Further, for $k\geq 2$ we have%
\begin{equation*}
I_{k}(t)=\frac{\left( 1+o(1)\right) }{\log ^{2k}t}\int_{2}^{t\varepsilon }%
\frac{\left( t-u\right) ^{k-1}}{\log u}du+\frac{\left( 1+o(1)\right) }{\log t%
}\int_{t\varepsilon }^{t-2}\frac{\left( t-u\right) ^{k-1}}{\log ^{2k}(t-u)}%
du.
\end{equation*}%
For large $t$%
\begin{equation*}
\int_{2}^{\varepsilon t}\frac{\left( t-u\right) ^{k-1}}{\log u}du\leq
2t^{k-1}\frac{\varepsilon t}{\log t}
\end{equation*}%
while%
\begin{equation*}
\int_{\varepsilon t}^{t-2}\frac{\left( t-u\right) ^{k-1}}{\log ^{2k}(t-u)}%
du=\int_{2}^{t(1-\varepsilon )}\frac{u^{k-1}}{\log ^{2k}u}du\sim \frac{1}{k}%
\frac{(1-\varepsilon )^{k}t^{k}}{\log ^{2k}t}.
\end{equation*}%
Hence letting first $t\rightarrow \infty $ and than $\varepsilon
\rightarrow +0$ the statement follows. $\Box$

\section{Proofs of the main results}\label{s:proofmain}

First we study the asymptotic behavior of the moments
\begin{equation*}
P_{k}(t):=\mathbf{E}\mu _{0}^{\left[ k\right] }(t),
\end{equation*}%
where $x^{\left[ k\right] }=x(x-1)\cdot \cdot \cdot (x-k+1),$ assuming that
the offspring generating function is infinitely differentiable at point $s=1$%
. To this aim we need the classical Faa di Bruno formula for the $n$-th
derivative of the composition of functions $g(h(s))$ (see \cite{FdB1857})$:$%
\begin{eqnarray*}
&&\frac{d^{n}g(h(s))}{ds^{n}} =g^{\prime}(h(s))h^{(n)}(s) \\
&&\quad +\sum_{k=2}^{n}g^{(k)}(h(s))\sum_{\substack{ j_{1}+j_{2}+\cdot \cdot
\cdot +j_{n-1}=k  \\ j_{1}+2j_{2}+\cdot \cdot \cdot +\left( n-1\right)
j_{n-1}=n}}\frac{n!}{j_{1}!j_{2}!\cdot \cdot \cdot j_{n-1}!}\left( \frac{%
h^{(1)}(s)}{1!}\right) ^{j_{1}}\cdot \cdot \cdot \left( \frac{h^{(n-1)}(s)}{%
\left( n-1\right) !}\right) ^{j_{n-1}}.
\end{eqnarray*}

\begin{lemma}
\label{Lmoment}If $f(s)$ satisfies Hypothesis (II) and is infinitely
differentiable at point $s=1$ then for any fixed $n\ge 1$,
\begin{equation}
P_{n}(t)\sim n!\left( \frac{\alpha f^{(2)}(1)}{2}\right) ^{n-1}\frac{1}{%
c_{4}^{2n-1}}\frac{t^{n-1}}{\log ^{2n-1}t}, \qquad t\to\infty.
\label{MomEst}
\end{equation}
\end{lemma}

\noindent\textbf{Proof.} From formula (\ref{ur20}) by
differentiation we have for all $t \ge 0$,
\begin{equation*}
P_{1}(t)=1-G_{1}(t)+\int_{0}^{t}\alpha f^{\prime}(1)P_{1}(t-u)\,dG_{1}(u)+\int_{0}^{t}(1-\alpha
)(1-h_{4})P_{1}(t-u)\,d(G_{1}\ast G_{2}(u)).
\end{equation*}

\noindent It follows that
\begin{equation*}
P_{1}(t)=1-G_{1}(t)+\int_{0}^{t}P_{1}(t-u)\,dK(u)=\left( 1-G_{1}\right) \ast
V_{K}(t)
\end{equation*}%
where
\begin{equation*}
V_{K}(t)=\sum_{j=0}^{\infty }K^{\ast j}(t)
\end{equation*}%
is the renewal function corresponding to the distribution function $K(t)$
(Recalling that $K(t):= \alpha f^{\prime}(1) G_1(t) + (1-\alpha) (1-h_4) G_1
\ast G_2(t)$). Since%
\begin{equation*}
1-K(t)\sim c_{4}t^{-1}
\end{equation*}%
we have by Theorem 3 in \cite{Eri70} that
\begin{equation}
P_{1}(t)\sim \frac{1}{c_{4}\log t}, \qquad t\to\infty.  \label{FirstMoment}
\end{equation}%
Note that in view of $G_{1}(t)=1-e^{-t}$ we have%
\begin{equation}
\frac{d}{dt}(G_{1}\ast V_{K}(t))=\left( 1-G_{1}\right) \ast
V_{K}(t)=P_{1}(t).  \label{DifV}
\end{equation}

Further, by writing  $F^{(n)}(t;s)= {\frac{\partial^n F(t; s)
}{\partial^n s}}$ for $n\geq 2$, we have
\begin{equation*}
F^{(n)}(t;s)=\int_{0}^{t}\alpha \frac{d^{n}f(F(t-u;s))}{ds^{n}}%
\,dG_{1}(u)+\int_{0}^{t}(1-\alpha )(1-h_{4})F^{(n)}(t-u;s)\,d(G_{1}\ast
G_{2}(u))
\end{equation*}
or, by the Faa di Bruno formula at $s=1$,
\begin{equation}
P_{n}(t)=\alpha
\int_{0}^{t}H_{n}(t-u)\,dG_{1}(u)+\int_{0}^{t}P_{n}(t-u)\,dK(u),  \label{Re2}
\end{equation}%
where%
\begin{equation}
H_{n}(t):=\sum_{k=2}^{n}f^{(k)}(1)\sum_{\substack{ j_{1}+j_{2}+\cdot \cdot
\cdot +j_{n-1}=k  \\ j_{1}+2j_{2}+\cdot \cdot \cdot +\left( n-1\right)
j_{n-1}=n}}\frac{n!}{j_{1}!j_{2}!\cdot \cdot \cdot j_{n-1}!}\left( \frac{%
P_{1}(t)}{1!}\right) ^{j_{1}}\cdot \cdot \cdot \left( \frac{P_{n-1}(t)}{%
\left( n-1\right) !}\right) ^{j_{n-1}}.  \label{DefH}
\end{equation}
Solving the renewal equation (\ref{Re2}) with respect to $P_{n}(t)$ gives
\begin{equation*}
P_{n}(t)=\alpha \int_{0}^{t}H_{n}(t-u)\,d\left( G_{1}\ast V_{K}(u)\right)
\end{equation*}%
or, in view of (\ref{DifV})
\begin{equation*}
P_{n}(t)=\alpha \int_{0}^{t}H_{n}(t-u)\,P_{1}(u)du.
\end{equation*}%
For $n=2$ we have%
\begin{equation*}
P_{2}(t)=\alpha f^{(2)}(1)\int_{0}^{t}P_{1}^{2}(t-u)\,P_{1}(u)du.
\end{equation*}%
On account of Lemma \ref{LbechInteg} this leads to%
\begin{equation*}
P_{2}(t)\sim \frac{\alpha f^{(2)}(1)}{c_{4}^{3}}I_{1}(t)\sim \frac{\alpha
f^{(2)}(1)}{c_{4}^{3}}\frac{t}{\log ^{3}t}=2!\frac{\alpha f^{(2)}(1)}{%
2c_{4}^{3}}\frac{t}{\log ^{3}t}.
\end{equation*}%
Now we use induction. Assume that for all $i<n$%
\begin{equation*}
P_{i}(t)\sim i!\frac{\left( \alpha f^{(2)}(1)\right) ^{i-1}}{%
2^{i-1}c_{4}^{2i-1}}\frac{t^{i-1}}{\log ^{2i-1}t}.
\end{equation*}%
Then for $2\leq k\leq n$ as $t\rightarrow \infty $
\begin{eqnarray*}
&&\sum_{\substack{ j_{1}+j_{2}+\cdot \cdot \cdot +j_{n-1}=k  \\ %
j_{1}+2j_{2}+\cdot \cdot \cdot +\left( n-1\right) j_{n-1}=n}}\frac{n!}{%
j_{1}!j_{2}!\cdot \cdot \cdot j_{n-1}!}\left( \frac{P_{1}(t)}{1!}\right)
^{j_{1}}\cdot \cdot \cdot \left( \frac{P_{n-1}(t)}{\left( n-1\right) !}%
\right) ^{j_{n-1}} \\
&\sim &\sum_{\substack{ j_{1}+j_{2}+\cdot \cdot \cdot +j_{n-1}=k  \\ %
j_{1}+2j_{2}+\cdot \cdot \cdot +\left( n-1\right) j_{n-1}=n}}\frac{n!}{%
j_{1}!j_{2}!\cdot \cdot \cdot j_{n-1}!}\left( \frac{1}{c_{4}\log t}\right)
^{j_{1}}\cdot \cdot \cdot \left( \frac{\left( \alpha f^{(2)}(1)\right) ^{n-2}%
}{2^{n-2}c_{4}^{2n-3}}\frac{t^{n-2}}{\log ^{2n-3}t}\right) ^{j_{n-1}} \\
&=&\left( \frac{\alpha f^{(2)}(1)}{2}\right) ^{n-k}\frac{1}{c_{4}^{2n-k}}%
\frac{t^{n-k}}{\log ^{2n-k}t}\sum_{\substack{ j_{1}+j_{2}+\cdot \cdot \cdot
+j_{n-1}=k  \\ j_{1}+2j_{2}+\cdot \cdot \cdot +\left( n-1\right) j_{n-1}=n}}%
\frac{n!}{j_{1}!j_{2}!\cdot \cdot \cdot j_{n-1}!}.
\end{eqnarray*}%
One may check that for any $n\ge 2$
\begin{equation*}
\sum_{\substack{ j_{1}+j_{2}+\cdot \cdot \cdot +j_{n-1}=2  \\ %
j_{1}+2j_{2}+\cdot \cdot \cdot +\left( n-1\right) j_{n-1}=n}}\frac{1}{%
j_{1}!j_{2}!\cdot \cdot \cdot j_{n-1}!}=\frac{n-1}{2}.
\end{equation*}
Thus, as $t\rightarrow \infty $%
\begin{eqnarray*}
H_{n}(t) &\sim &f^{(2)}(1)\left( \frac{\alpha f^{(2)}(1)}{2}\right) ^{n-2}%
\frac{1}{c_{4}^{2n-2}}\frac{t^{n-2}}{\log ^{2n-2}t}\sum_{\substack{ %
j_{1}+j_{2}+\cdot \cdot \cdot +j_{n-1}=2  \\ j_{1}+2j_{2}+\cdot \cdot \cdot
+\left( n-1\right) j_{n-1}=n}}\frac{n!}{j_{1}!j_{2}!\cdot \cdot \cdot
j_{n-1}!} \\
&=&n!\frac{f^{(2)}(1)}{2}\left( \frac{\alpha f^{(2)}(1)}{2}\right) ^{n-2}%
\frac{(n-1)}{c_{4}^{2n-2}}\frac{t^{n-2}}{\log ^{2n-2}t}.
\end{eqnarray*}%
Therefore, on account of Lemma \ref{LbechInteg}%
\begin{eqnarray*}
P_{n}(t) &=&\alpha \int_{0}^{t}H_{n}(t-u)\,P_{1}(u)du \\
&\sim &n!\left( \frac{\alpha f^{(2)}(1)}{2}\right) ^{n-1}\frac{(n-1)}{%
c_{4}^{2n-1}}\int_{2}^{t-2}\frac{\left( t-u\right) ^{n-2}}{\log
^{2n-2}\left( t-u\right) }\frac{1}{\log u}du \\
&\sim &n!\left( \frac{\alpha f^{(2)}(1)}{2}\right) ^{n-1}\frac{1}{%
c_{4}^{2n-1}}\frac{t^{n-1}}{\log ^{2n-1}t},
\end{eqnarray*}%
as desired. $\Box$

\begin{corollary}
\label{Cor1}If $f(s)$ satisfies Hypothesis (II) then%
\begin{equation}
\lim \inf_{t\rightarrow \infty }\frac{tq(t)}{\log t}\geq \frac{c_{4}}{\alpha
f^{(2)}(1)}\text{.}  \label{qbelow}
\end{equation}
\end{corollary}

\noindent\textbf{Proof}.   From the proof of Lemma \ref{Lmoment} it
is clear that for the asymptotic representation (\ref{MomEst}) be
valid for $n=1,2$ it suffices that $f^{(2) }(1)<\infty $.
From this and the Lyapunov inequality
\begin{equation*}
q(t)=\mathbf{P}\left( \mu _{0}(t)>0\right) \geq \frac{\left( \mathbf{E}\mu
_{0}(t)\right) ^{2}}{\mathbf{E}\mu _{0}^{2}(t)}\sim \frac{\log t}{t}\frac{%
c_{4}}{\alpha f^{(2)}(1)}
\end{equation*}%
the needed statement easily follows. $\Box$

\medskip
Before giving the exact asymptotic of $q(t)$,    we show  at first a
rough upper bound for $q(t)$.

\begin{lemma}
\label{Lrough}We have%
\begin{equation*}
\lim \sup_{t\rightarrow \infty }\frac{tq(t)}{\log t}<\infty .
\end{equation*}
\end{lemma}

\noindent\textbf{Proof.} Fix an arbitrary $u>0$. Define for $0 \le x
\le 1$,
\begin{equation*}
T(x):=xk_{4}(u)-\alpha \Phi (x)e^{-u}=\alpha (f^{\prime}(1)x-\Phi
(x))e^{-u}+x(1-\alpha )(1-h_{4})(G_{1}\ast G_{2})^{\prime }(u).
\end{equation*}%

 Recalling  (\ref{AsymPhi}). We have
\begin{eqnarray*}
T^{(1)}(x) &=&\alpha  f^{\prime}(1-x) e^{-u}+(1-\alpha
)(1-h_{4})(G_{1}\ast G_{2})^{\prime }(u)   >0, \quad 0\le x \le 1.
\end{eqnarray*}%
  Hence
 $T(x)$ is monotone increasing in $x\in (0,1)$.   This fact will be used several times in the sequel.

Let us write a formal
representation
\begin{equation*}
q(t)=\beta (t)\frac{\log \left( t+1\right) }{t+1}, \qquad t\ge 0,
\end{equation*}%
and set
\begin{equation*}
q(u,t):=\beta (t)\frac{\log \left( t-u+1\right) }{t-u+1}, \qquad 0 \le u \le t.
\end{equation*}%
Assume that the desired statement is not true, that is that
\begin{equation}
\lim \sup_{t\rightarrow \infty } \beta(t)= \infty.  \label{limsupbeta}
\end{equation}%

Under (\ref{limsupbeta}), there exists a sequence $t_{k}\rightarrow
\infty $ as $k\rightarrow \infty $   such that $\beta
(t_{k})=\sup_{u\leq t_{k}}\beta (u)$ and $\lim_{k\to\infty}
\beta(t_k)= \infty$. Then, for sufficiently large $t=t_{k}$ (we omit
the low index for simplicity), $\beta(t)=\sup_{0\le s\le t} \beta(s)
>1$. In view of (\ref{LintStar}) with $p=1$, we have
\begin{eqnarray}
&&\left\vert \int_{t/2}^{t}q(t-u)\left( k_{4}(u)-\alpha \Psi
(q(t-u))e^{-u}\right) du\right\vert   \nonumber\\
&\le &\int_{t/2}^{t}\beta (t-u)\frac{\log \left( t-u+1\right) }{t-u+1}%
k_{4}(u)du+\alpha \int_{t/2}^{t}\Phi (q(t-u))e^{-u}du  \nonumber\\
&\leq &\beta (t)\frac{ 2 c_{4}}{t^{2}}\log ^{2}t,  \label{t/2}
\end{eqnarray}%

 \noindent since $\beta(t)>1$ and $\alpha \int_{t/2}^{t}\Phi (q(t-u))e^{-u}du \le \alpha f'(1) e^{ - t/2}= o(\frac{\log ^2 t}{t^2})$.

 Let us observe that  $q(t) \to 0$ as $t\to\infty$. In fact, we have a rough bound:   $q(t) = {\bf P}(\mu_0(t)>0) \le   {\bf E} \mu_0(t)= P_1(t) \sim {1 \over c_4 \log t}$ by (\ref{FirstMoment}). It
follows that for  $0\leq u\leq t/2$,
\begin{eqnarray*}
q(t-u) &=&\beta (t-u)\frac{\log \left( t-u+1\right) }{t-u+1}\leq \beta (t)%
\frac{\log \left( t-u+1\right) }{t-u+1} \\
&=&q(u,t)\leq \beta (t)\frac{\log \left( t/2\right) }{t/2} \le 3  \beta (t)\frac{\log t}{t}= 3q(t) \leq 1,
\end{eqnarray*}%
for all sufficiently large $t$. Thus, the inequality
\begin{equation*}
q(t-u)\left( k_{4}(u)-\alpha \Psi (q(t-u))e^{-u}\right) \leq q(u,t)\left(
k_{4}(u)-\alpha \Psi (q(u,t))e^{-u}\right)
\end{equation*}%
 is valid for   $0\leq u\leq t/2$. This, in view of (\ref{Eqvq}) and (\ref{t/2}),  implies that
\begin{eqnarray*}
&& \beta (t)\frac{\log (t+1)}{t+1} =q(t) \\ &=&
1-G_{1}(t)+\int_{0}^{t}q(t-u)\left( k_{4}(u)-\alpha \Psi
(q(t-u))e^{-u}\right) du   \\
&\le &\beta (t)\frac{ 3 c_{4}}{t^{2}}\log ^{2}t \\
&&+\int_{0}^{t/2}\beta (t-u)\frac{\log \left( t-u+1\right) }{t-u+1}%
\left( k_{4}(u)-\alpha \Psi \left( \beta (t-u)\frac{\log \left( t-u+1\right)
}{t-u+1}\right) e^{-u}\right) du \\
&\leq &\beta (t)\frac{ 3 c_{4}}{t^{2}}\log ^{2}t  \\
&&+\beta (t)\int_{0}^{t/2}\frac{\log \left( t-u+1\right) }{t-u+1}%
\left( k_{4}(u)-\alpha \Psi \left( \beta (t)\frac{\log \left( t-u+1\right) }{%
t-u+1}\right) e^{-u}\right) du.
\end{eqnarray*}%
By (\ref{Lintegr}), we have that
\begin{equation}
\int_{0}^{t/2}\frac{\log \left( t-u+1\right) }{t-u+1}k_{4}(u)du \le \frac{%
\log (t+1)}{t+1} +  2 c_4\, \frac{\log ^{2}t}{t^{2}}   \label{IntK}
\end{equation}%
and, for any small  $\delta >0$ and sufficiently large $t$,
\begin{eqnarray*}
&&\int_{0}^{t/2}\frac{\log \left( t-u+1\right) }{t-u+1}\alpha \Psi
\left( \beta (t)\frac{\log \left( t-u+1\right) }{t-u+1}\right) e^{-u}du \\
&\geq &(1-\delta )\beta (t)\frac{\alpha f^{(2)}(1)}{2}\int_{0}^{t/2}%
\frac{\log ^{2}\left( t-u+1\right) }{\left( t-u+1\right) ^{2}}e^{-u}du\geq
(1-2\delta )\beta (t)\frac{\alpha f^{(2)}(1)}{2}\frac{\log ^{2}t}{t^{2}}.
\end{eqnarray*}
As a result we get that
$$
\beta (t)\frac{\log (t+1)}{t+1} \leq \beta (t)\left( \frac{\log (t+1)}{t+1}%
+ 5 c_4\frac{\log ^{2}t}{t^{2}}\right)  -(1-2\delta )\beta ^{2}(t)\frac{\alpha f^{(2)}(1)}{2}\frac{\log ^{2}t}{%
t^{2}}  .$$

\noindent
Hence, after simplification we see that%
\begin{equation*}
(1-2\delta )\beta ^{2}(t)\frac{\alpha f^{(2)}(1)}{2}\leq  5 c_4 \beta (t),
\end{equation*}%
which  is impossible if $\beta (t)\rightarrow \infty $. Hence $ \limsup_{t \to\infty} \beta(t) < \infty$ and the lemma is proved. $\Box$

The crucial step in the proof of   Theorem 1 is the following lemma:

\begin{proposition}
\label{Lprecise}In the case $d=4$,
\begin{equation}
q(t)=C\frac{\log t}{t}(1+o(1)),  \label{Asympq}
\end{equation}%
where
\begin{equation}
C:=\frac{3c_{4}}{\alpha f^{(2)}(1)}.\label{ddefC}
\end{equation}
\end{proposition}

\noindent\textbf{Proof.} We will use a formal representation%
\begin{equation*}
q(t)=\frac{\log (t+1)}{t+1}\left( C+\frac{\rho (t)}{\sqrt{\log (t+1)}}%
\right) =C\frac{\log (t+1)}{t+1}+\frac{\rho (t)\sqrt{\log (t+1)}}{t+1}.
\end{equation*}%
Note that by Corollary \ref{Cor1}
\begin{equation}
\lim \inf_{t\rightarrow \infty }\frac{\rho (t)}{\sqrt{\log t}}=C_{0}\geq -%
\frac{2}{3}C,  \label{RhoBelow}
\end{equation}%
and, by Lemma \ref{Lrough}%
\begin{equation*}
\lim \sup_{t\rightarrow \infty }\frac{\rho (t)}{\sqrt{\log t}}\leq
C_{1}<\infty ,
\end{equation*}
for some constant $C_1>0$.  If $\rho (t)$ is bounded then the lemma
is proved. Thus, assume that
\begin{equation*}
\lim \sup_{t\rightarrow \infty }\rho (t)=\infty \text{ and }\lim
\inf_{t\rightarrow \infty }\rho (t)=-\infty .\text{ }
\end{equation*}%
Clearly, if we prove the desired statement for this case then the cases when
one of the limits above is finite will follow easily.

Assume first that $\lim \sup_{t\rightarrow \infty }\rho (t)=\infty $ and let
\begin{equation*}
\mathcal{A}_{+}:=\left\{ t\in \lbrack 0,\infty ):\rho (t)=\sup_{v\leq t}\rho
(v)\right\} .
\end{equation*}

\noindent Then there exists an unbounded  sequence $(t_k) \in {\cal
A}_+$ such that $\rho(t_k) \to \infty$ as $k \to \infty$. Our
subsequent arguments are for large  $t\in \mathcal{A}_{+}$. Then a
priori, $\rho(t)>1$. For
 $0\leq u\leq t$,   set
\begin{equation*}
Q(u,t):=C\frac{\log (t-u+1)}{t-u+1}+\frac{\rho (t)\sqrt{\log (t-u+1)}}{%
t-u+1}\geq q(t-u).
\end{equation*}%

Fix a small  $\varepsilon>0$.  Clearly, for $0\leq u\leq
t\varepsilon $ and sufficiently large $t$,
\begin{equation*}
Q(u,t)\leq C\frac{\log (t+1)}{t(1-\varepsilon )}+\frac{\rho
(t)\sqrt{\log \left( t+1\right) }}{t(1-\varepsilon )}\leq
\frac{1}{1-\varepsilon }q(t)< 1,
\end{equation*}

\noindent since $q(t) \to 0$. Using again the monotonicity of the
function: $x (\in (0,1)) \to xk_{4}(u)-\alpha \Phi (x)e^{-u}$ (just
like in Lemma \ref{Lrough}),  we get the inequality
\begin{equation*}
q(t-u)\left( k_{4}(u)-\alpha \Psi \left( q(t-u)\right) e^{-u}\right)
\leq Q(u,t)\left( k_{4}(u)-\alpha \Psi \left( Q(u,t)\right)
e^{-u}\right).
\end{equation*}

\noindent It follows that
\begin{eqnarray}
q(t) &=&\int_{t\varepsilon }^{t}q(t-u)k_{4}(u)du  \notag \\
&&+\int_{0}^{t\varepsilon }q(t-u)\left( k_{4}(u)-\alpha \Psi \left(
q(t-u)\right) e^{-u}\right) du+ o\left(\frac{\log^2 t}{t^{2}}\right)  \notag \\
&\leq &\int_{t\varepsilon
}^{t}Q(u,t)k_{4}(u)du+\int_{0}^{t\varepsilon }Q(u,t)\left(
k_{4}(u)-\alpha \Psi \left( Q(u,t)\right) e^{-u}\right) du +
o\left(\frac{\log^2 t}{t^{2}}\right)
\notag \\
&=&\int_{0}^{t}Q(u,t)k_{4}(u)du-\alpha \int_{0}^{t\varepsilon }\Phi
\left( Q(u,t)\right) e^{-u}du+ o\left(\frac{\log^2 t}{t^{2}}\right).
\label{Basic1}
\end{eqnarray}

We now evaluate the integrals in (\ref{Basic1}). Recalling
(\ref{CorInt}) with $p=1$ and $p=1/2$, we have
\begin{eqnarray*}
&&\int_{0}^{t}Q(u,t)k_{4}(u)du =C\int_{0}^{t}\frac{\log (t-u+1)}{t-u+1}
k_{4}(u)du+\rho (t)\int_{0}^{t}\frac{\sqrt{\log \left( t-u+1\right)
}
}{t-u+1} k_{4}(u)du \\
&&\qquad=C\frac{\log (t+1)}{t+1}+\frac{3Cc_{4}\log ^{2}t}{2t^{2}}(1+o(1)) +\frac{\rho (t)\sqrt{\log (t+1)}}{t+1}+\frac{5c_{4}\rho (t)\log ^{3/2}t}{%
3t^{2}}(1+o(1)).
\end{eqnarray*}
Since all the moments of the distribution with density $e^{-u}$ are finite,
and $u^{k}e^{-u}$ decays rapidly at infinity for any $k\geq 0,$ we have
\begin{eqnarray}
&&-\alpha \int_{0}^{t\varepsilon }\Phi \left( Q(u,t)\right) e^{-u}du=-\alpha \Phi \left( q(t)\right) (1+o(1))=-\alpha \frac{f^{(2)}(1)}{2}q^{2}(t)(1+o(1)) \nonumber\\
&&\qquad\qquad=-\alpha \frac{f^{(2)}(1)}{2}\left( C^{2}\frac{\log ^{2}t}{t^{2}}+2C\frac{%
\rho (t)\log ^{3/2}t}{t^{2}}+\frac{\rho ^{2}(t)\log t}{t^{2}}\right)
(1+o(1)) \nonumber
\\
&&\qquad\qquad=-\alpha \frac{f^{(2)}(1)}{2}\left( C^{2}\frac{\log ^{2}t}{t^{2}}+2C\frac{%
\rho (t)\log ^{3/2}t}{t^{2}}+\frac{\rho ^{2}(t)\log t}{t^{2}}\right)
+o\left( \frac{\log ^{2}t}{t^{2}}\right)  \label{phiqut}.
\end{eqnarray}
Substituting this in (\ref{Basic1}) we get that for any small
$\delta>0$,
\begin{eqnarray*}
&& C\frac{\log (t+1)}{t+1}+\frac{\rho (t)\sqrt{\log (t+1)}}{t+1}
 = q(t) \\ &&\qquad\leq \delta
\, \frac{\log ^{2}t}{t^{2}}+C\frac{\log (t+1)}{t+1}+\frac{3Cc_{4}\log ^{2}t}{%
2t^{2}} +\frac{\rho (t)\sqrt{\log (t+1)}}{t+1}+\frac{5c_{4}\rho (t)\log ^{3/2}t}{%
3t^{2}} \\
&&\qquad\qquad-\alpha \frac{f^{(2)}(1)}{2}\left( C^{2}\frac{\log ^{2}t}{t^{2}}+2C\frac{%
\rho (t)\log ^{3/2}t}{t^{2}}+\frac{\rho ^{2}(t)\log t}{t^{2}}\right)
\end{eqnarray*}%
which, on account the definition (\ref{ddefC})
leads, after natural transformations, to

\begin{equation*}
\frac{\alpha f^{(2)}(1 ) \rho (t)}{2\sqrt{\log t}}\left( \frac{8}{9}C+\frac{%
\rho (t)}{\sqrt{\log t}}\right) \leq \delta.
\end{equation*}

Recall that $\rho(t)>1$. We get that $\frac{\alpha f^{(2)}(1 ) \rho
(t)}{2\sqrt{\log t}}    \leq  {9\delta \over 8}$ and hence
\begin{equation*}
\lim \sup_{t\rightarrow \infty ,t\in \mathcal{A}_{+}}\frac{\rho (t)}{\sqrt{%
\log t}}=0.
\end{equation*}%
Let
\begin{equation*}
t_{+}(u)=\sup \left\{ t\leq u:t\in \mathcal{A}_{+}\right\} .
\end{equation*}%
Clearly,%
\begin{equation*}
\rho (u)\leq \rho (t_{+}(u)).
\end{equation*}%
Now%
\begin{eqnarray}
\lim \sup_{u\rightarrow \infty }\frac{u q(u)}{\log u} &=&C+\lim
\sup_{u\rightarrow \infty }\left( \frac{\rho (u)}{\sqrt{\log
u}}\right) \leq C+\lim \sup_{u\rightarrow \infty }\left( \frac{\rho
(t_{+}(u))}{\sqrt{\log
t_{+}(u)}}\right)   \notag \\
&=&C+\lim \sup_{t\rightarrow \infty ,t\in \mathcal{A}_{+}}\left( \frac{\rho
(t)}{\sqrt{\log t}}\right) =C.  \label{RR1}
\end{eqnarray}

To get estimate from below we assume that%
\begin{equation*}
\lim \inf_{t\rightarrow \infty }\rho (t)=-\infty
\end{equation*}
(otherwise, we are done) and let
\begin{equation*}
\mathcal{A}_{-}:=\left\{ t\in \lbrack 0,\infty ):\rho (t)=\inf_{u\leq t}\rho
(v)\right\} .
\end{equation*}
Our subsequent arguments are for large  $t\in \mathcal{A}_{-}$. For
$0\leq u\leq t$
\begin{equation*}
Q(u,t)=C\frac{\log (t-u+1)}{t-u+1}+\frac{\rho (t)\sqrt{\log (t-u+1)}}{t-u+1}%
\leq q(t-u)\leq 1.
\end{equation*}

Fix a small $\varepsilon>0$ and consider   $0\leq u\leq
t\varepsilon$.  In view of (\ref{RhoBelow}), we get that
\begin{equation*}
Q(u,t)\geq C\frac{\log (t+1)}{t+1}+\frac{\rho (t)\sqrt{\log t(1-\varepsilon )%
}}{t(1-\varepsilon )}\geq \frac{1}{4}C\frac{\log
(t+1)}{t+1}>0\text{.}
\end{equation*}

\noindent  Just like in (\ref{Basic1}), we use again the
monotonicity in the reverse order and get  that
$$
q(t) \ge \int_{0}^{t}Q(u,t)k_{4}(u)du-\alpha \int_{0}^{t\varepsilon
}\Phi \left( Q(u,t)\right) e^{-u}du+o\Big(\frac{\log^2
t}{t^{2}}\Big).$$

By (\ref{CorInt}) with $p=1$ and $p=1/2$ (and using the fact that
$|\rho(t)|= O( \sqrt{\log t})$),  we get as before
\begin{eqnarray*}
 \int_{0}^{t}Q(u,t)k_{4}(u)du&=&C\frac{\log (t+1)}{t+1}+\frac{3Cc_{4}\log ^{2}t}{2t^{2}}
  \\&&+\frac{\rho (t)%
\sqrt{\log (t+1)}}{t+1}+\frac{5c_{4}\rho (t)\log ^{3/2}t}{3t^{2}}
 +o\left( \frac{\log ^{2}t}{t^{2}}\right) .
\end{eqnarray*}%

\noindent This together with (\ref{phiqut}) implies that  for any
small $\delta>0$,
\begin{eqnarray*}
&& C\frac{\log (t+1)}{t+1}+\frac{\rho (t)\sqrt{\log (t+1)}}{t+1} = q(t) \\
 &&\qquad\geq C\frac{%
\log (t+1)}{t+1}+\frac{3Cc_{4}\log ^{2}t}{2t^{2}}  +\frac{\rho (t)\sqrt{\log (t+1)}}{t+1}+\frac{5c_{4}\rho (t)\log ^{3/2}t}{%
3t^{2}} \\
&&\qquad\quad-\alpha \frac{f^{(2)}(1)}{2}\left( C^{2}\frac{\log ^{2}t}{t^{2}}+2C\frac{%
\rho (t)\log ^{3/2}t}{t^{2}}+\frac{\rho ^{2}(t)\log t}{t^{2}}\right)
- \delta  \frac{\log ^{2}t}{t^{2}}
\end{eqnarray*}%
which, similarly to the previous case gives after simplifications
\begin{equation*}
\delta \geq -\frac{\rho (t)\alpha f^{(2)}(1)}{2\sqrt{\log t}}\left( \frac{8}{9}C+%
\frac{\rho (t)}{\sqrt{\log t}}\right)  .
\end{equation*}%

This,   on account of (\ref{RhoBelow}) gives that ($t$ being
  large)
\begin{equation*}
\delta \geq -\frac{\rho (t)\alpha f^{(2)}(1)}{2\sqrt{\log
t}}\frac{1}{9}C .
\end{equation*}

\noindent Since $\delta$ is arbitrary, we get that
\begin{equation*}
\lim \sup_{t\rightarrow \infty ,t\in \mathcal{A}_{-}}\frac{\left\vert \rho
(t)\right\vert }{\sqrt{\log t}}=0.
\end{equation*}%
Let
\begin{equation*}
t_{-}(u)=\sup \left\{ t\leq u:t\in \mathcal{A}_{-}\right\} .
\end{equation*}%
Clearly,%
\begin{equation*}
\left\vert \rho (u)\right\vert \leq \left\vert \rho (t_{-}(u))\right\vert .
\end{equation*}%
Now%
\begin{eqnarray}
\lim \inf_{u\rightarrow \infty }\frac{uq(u)}{\sqrt{\log u}} &=&C-\lim
\sup_{u\rightarrow \infty }\left( \frac{\left\vert \rho (u)\right\vert }{%
\sqrt{\log u}}\right) \geq C-\lim \sup_{u\rightarrow \infty }\left( \frac{%
\left\vert \rho (t_{-}(u))\right\vert }{\sqrt{\log t_{-}(u)}}\right)  \notag
\\
&=&C-\lim \sup_{t\rightarrow \infty ,t\in \mathcal{A}_{-}}\frac{\left\vert
\rho (t)\right\vert }{\sqrt{\log t}}=C.  \label{RR2}
\end{eqnarray}

Combining (\ref{RR1}) and (\ref{RR2}) \ gives
\begin{equation*}
\lim_{t\rightarrow \infty }\frac{tq(t)}{\log t}=C,
\end{equation*}%
proving the Proposition. $\Box$

\begin{theorem}
\label{Tonedim4}Assume that $f(s)$ is infinitely differentiable at point $%
s=1 $ and satisfies%
\begin{equation*}
\alpha f^{\prime }(1)+(1-\alpha )(1-h_{4})=1.
\end{equation*}%
Then%
\begin{equation*}
\lim_{t\rightarrow \infty }\mathbf{P}\left( \frac{\mu _{0}(t)q(t)}{P_{1}(t)%
}\leq x\left\vert \mu _{0}(t)>0\right. \right) =\frac{1}{3}+\frac{2}{3}%
\left( 1-e^{-2x/3}\right) ,x>0,
\end{equation*}%
or, what is the same%
\begin{equation*}
\lim_{t\rightarrow \infty }\mathbf{E}\left[ e^{-\lambda \mu
_{0}(t)q(t)/P_{1}(t)}|\mu _{0}(t)>0\right] =\frac{1}{3}+%
\frac{2}{3}\frac{2}{2+3\lambda },\lambda \geq 0.
\end{equation*}
\end{theorem}

\noindent\textbf{Proof}. It follows from Lemma \ref{Lmoment} that%
\begin{eqnarray*}
\mathbf{E}\mu _{0}^{n}(t) &\sim &n!\left( \frac{\alpha f^{(2)}(1)}{2}\right)
^{n-1}\frac{1}{c_{4}^{2n-1}}\frac{t^{n-1}}{\log ^{2n-1}t} \\
&=&\left( \frac{3}{2}\right) ^{n-1}n!\left( \frac{\alpha f^{(2)}(1)}{3c_{4}}%
\frac{t}{\log t}\right) ^{n-1}\left( \frac{1}{c_{4}\log t}\right) ^{n} \\
&\sim &\left( \frac{3}{2}\right) ^{n-1}n!\frac{P_{1}^{n}(t)}{q^{n-1}(t)}=%
q(t)\left( \frac{3}{2}\right) ^{n-1}n!\left( \frac{P_{1}(t)}{q(t)}\right) ^{n}.
\end{eqnarray*}%
Therefore, as $t\rightarrow \infty $%
\begin{equation}
\mathbf{E}\left[ \left( \frac{\mu _{0}(t)q(t)}{P_{1}(t)}\right)
^{n}\left\vert \mu _{0}(t)>0\right. \right] =\frac{1}{q(t)}\left( \frac{q(t)%
}{P_{1}(t)}\right) ^{n}\mathbf{E}\mu _{0}^{n}(t)\rightarrow \frac{2}{3}\left( \frac{3}{2}\right) ^{n}n!.
\label{MomConv}
\end{equation}%
Thus, for any $n\geq 1$ the $n$-th moment of the conditional
distribution converges to the $n-$th moment of the mixture (with
probabilities 2/3 and 1/3, respectively) of the exponential
distribution with parameter $2/3$ (which is uniquely defined by its
moments) and the distribution having the unit atom at zero. Hence
the statement of the theorem follows. $\Box$

\medskip
Our next step is to generalize Theorem \ref{Tonedim4} to the case of
arbitrary probability generating function $f(s)$ with finite second moment.
To this aim we need an approximation lemma.

\begin{lemma}
\label{Lapprox} Let $f(s)$ be an arbitrary probability generating function
with $f^{\prime }(1)>0$ and $f^{(2) }(1)\in (0,\infty )$. For any $%
\varepsilon \in (0,1)$, there exist two polynomial probability generating
functions $f_{-}(s)$, $f_{+}(s)$ and some constant $s_{0}=s_{0}(f_{-},f_{+},%
\varepsilon )<1$ such that
\begin{equation}
f_{-}(s)\leq f(s)\leq f_{+}(s)  , \qquad \forall s\in (s_{0},1],
\label{IneqFunct}
\end{equation}%
 and
\begin{equation*}
f_{-}^{\prime }(1)=f_{+}^{\prime }(1)=f^{\prime }(1),
\end{equation*}%
and%
\begin{equation}
{\frac{f_{+}^{{(2)   }}(1)}{1+\varepsilon }}\leq f^{ (2) }(1)\leq {\frac{%
f_{-}^{{(2) }}(1)}{1-\varepsilon }}.  \label{DoublIneq}
\end{equation}
\end{lemma}

Remark that necessarily,  $ f_+^{{(2)}}(1) \ge f^{{(2) }}(1) \ge
f_-^{{(2) }}(1)$.

\noindent\textbf{Proof}. Let $N$ be an integer valued random
variable with generating function $f$. We only need to consider the
unbounded $N$ case, otherwise there is nothing to prove.

Let $\varepsilon>0$ be small. Assume for the moment that there exist two
integer valued and bounded random variables $N_1$ and $N_2$, such that
\begin{equation}  \label{en1}
\mathbf{E} (N_1)= \mathbf{E} (N_2)= \mathbf{E} (N), \qquad {\mbox var}(N_1)
< {\mbox var}(N) < {\mbox var}(N_2),
\end{equation}
and
\begin{equation}  \label{en2}
{\mbox var}(N_2) - \varepsilon < {\mbox var}(N ) < {\mbox var}(N_1) +
\varepsilon.
\end{equation}

\noindent Define $f_-(s):= \mathbf{E}(s^{N_1})$, $f_+(s):= \mathbf{\ E}%
(s^{N_2})$ for $0\le s\le 1$. Then (\ref{IneqFunct}) follows from (\ref{en1}%
) by  developing the three generating functions at $1$, whereas (\ref%
{DoublIneq}) follows from (\ref{en2})  since $\varepsilon$ is arbitrary.

To construct $N_1$ and $N_2$ satisfying (\ref{en1}) and (\ref{en2}), we fix $%
k $ an integer sufficiently large such that $0< \mathbf{E}( N (N-k)^+) \le
\varepsilon/3 $ and $\mathbf{E}((N-k)^-) \ge 1$ (where $x^+:=\max(x,0)$ and $x^-:=\max(-x, 0)$ for any real $x$).  Since $N= N \wedge k + (N-k)^+$,  elementary
computation
 shows that \begin{equation*}
{\mbox var}(N)={\mbox var}(N\wedge k)+   {\mbox var}((N-k)^{+})+ 2 \mathbf{E}
((N-k)^{+})\,\mathbf{E}((N-k)^{-}).
\end{equation*}

 \noindent Therefore,
\begin{equation}
2\mathbf{E}((N-k)^+) \le {\mbox var}(N )- {\mbox var}(N\wedge k) \le  2 \mathbf{
E}(N(N-k)^+) \le {\frac{2\varepsilon}{3}}. \label{varnk}
\end{equation}

Let $r=\mathbf{E}((N-k)^{+})>0$. Plainly $r \le {1\over k} \mathbf{E}( N (N-k)^+) \le
{\varepsilon\over 3k} < 1$.  Choose a Bernoulli variable $B$ with $%
\mathbf{P}(B=1)=r=1-\mathbf{P}(B=0)$, independent of $N$. Define $%
N_{1}:=N\wedge k+B$. Hence $\mathbf{E}(N_{1})=\mathbf{E}(N)$. On the other hand, it follows from (\ref{varnk}) that   ${\mbox var}(N )\le {\mbox var}(N\wedge k) +
2\varepsilon/3< {\mbox var}(N_1) + 2\varepsilon/3$, and ${\mbox
var}(N_1)= {\mbox var}(N\wedge k) + r( 1-r) \le {\mbox var}(N\wedge k) + r
\le {\mbox var}(N ) - r $ since $\mathbf{E}((N-k)^+)=r$. Then $N_1$ fulfills
the conditions in (\ref{en1}) and (\ref{en2}).

To construct $N_2$, we choose $\ell:= \lfloor 3{\frac{ \mathbf{E}( N
(N-k)^+) ) }{\mathbf{E}(  (N-k)^+)}} \rfloor$ and $b:= {\frac{1}{\ell}}
\mathbf{E}( (N-k)^+)$. Let $\widetilde B$ be a Bernoulli variable with $%
\mathbf{P}(\widetilde B=1)= b = 1- \mathbf{P}(\widetilde B=0)$, independent
of $N$. Define  $N_2:= N\wedge k + \ell\, \widetilde B$. Plainly, $\mathbf{E}%
(N_2)= \mathbf{\ E}(N) $ and
\begin{equation*}
{\mbox var}(N_2)= {\mbox var}(N\wedge k) + \ell^2  b (1-b) < {\mbox var}
(N\wedge k) + 3 \mathbf{E}( N (N-k)^+) ) \le {\mbox var}(N ) +\varepsilon.
\end{equation*}

\noindent Note that $\ell^2  b (1-b) \ge \frac{5}{2}\, \ell b
\,{\frac{ \mathbf{E}( N (N-k)^+) ) }{\mathbf{E}(  (N-k)^+)}} =
\frac{5}{2}\,\, \mathbf{E}( N (N-k)^+) ) $. It follows that  ${\mbox
var}(N_2) \ge {\mbox var}(N \wedge k) + \frac{5}{2}\,\, \mathbf{E} (
N (N-k)^+) ) \ge {\mbox var}(N )+ \frac{1}{2}\,\mathbf{E}( N
(N-k)^+) ) > {\mbox var}(N ).$ This shows that  $N_2$ also fulfills
the conditions in (\ref{en1}) and (\ref{en2}) and completes the
proof of lemma. $\Box$

\begin{lemma}
\label{Lccomp}If there are two branching random walks with branching at the origin only whose
offspring generating functions $f_{1}(s)$ and $f_{2}(s)$ are such that%
\begin{equation*}
\alpha f_{i}^{\prime }(1)+(1-\alpha )(1-h_{4})=1,\quad i=1,2,
\end{equation*}%
and for some constant $0\le s_0<1$,
\begin{equation*}
f_{1}(s)\leq f_{2}(s), \qquad \forall\, s_0<s \le1,
\end{equation*}%
 then the respective generating functions $F^{1}(t;s)$
and $F^{2}(t;s)$ for the number of particles at the origin at moment $t$
meet the inequality%
\begin{equation}
F^{1}(t;s)\leq F^{2}(t;s)  \label{CompFF}
\end{equation}%
for all $s\in (s_{0},1]$.
\end{lemma}

\noindent\textbf{Proof}.  Let  $0\le s\le1$ and $t\ge0$. Introduce the notation%
\begin{eqnarray}
L\left( f,F\right) ( t;s) &:=&s(1-G_{1}(t))+(1-\alpha
)(1-h_{4})(1-G_{2}(\cdot ))\ast G_{1}(t)  \notag \\
&+&\int_{0}^{t}\alpha f(F(t-u;s))\,dG_{1}(u)+(1-\alpha )h_{4}G_{1}(t)  \notag
\\
&+&\int_{0}^{t}(1-\alpha )(1-h_{4})F(t-u;s)\,d(G_{1}\ast G_{2}(u)),
\end{eqnarray}

\noindent and for $i=1,2$,  set
\begin{equation*}
F_{0}^{i}(t;s)=s, \quad F_{n+1}^{i}(t;s):=L\left(
f_{i},F_{n}^{i}\right) \left( t;s\right) .
\end{equation*}

Let us show by induction on $n$  that
\begin{equation*}
s\leq F_{n}^{i}(t;s)\leq F_{n+1}^{i}(t;s), \qquad \forall 0\le s\le
1.
\end{equation*}

 Indeed, the expression%
\begin{equation*}
R_{i}(s):=\alpha f_{i}(s)+(1-\alpha )(1-h_{4})s+(1-\alpha )h_{4}
\end{equation*}%
is a probability generating function with $R_{i}^{\prime }(1)=\alpha
f_{i}^{\prime }(1)+(1-\alpha )(1-h_{4})=1$. Hence $R_{i}(s)\geq s$ for all $%
s\in \lbrack 0,1]$. Using this fact, we have%
\begin{eqnarray*}
&& F_{1}^{i}(t;s) =L\left( f_{i},F_{0}^{i}\right)  ( t;s )
 \\ &=&s(1-G_{1}(t))+(1-\alpha )(1-h_{4})(1-G_{2}(\cdot ))\ast G_{1}(t) \\
&&+\int_{0}^{t}\alpha f_{i}(s)\,dG_{1}(u)+(1-\alpha )h_{4}G_{1}(t)+\int_{0}^{t}(1-\alpha )(1-h_{4})s\, d(G_{1}\ast G_{2}(u)) \\
&=&s(1-G_{1}(t))+(1-\alpha )(1-h_{4})(1-G_{2}(\cdot ))\ast
G_{1}(t)+(1-\alpha )h_{4}G_{1}(t) \\
&&+\alpha f_{i}(s)G_{1}(t)+(1-\alpha )(1-h_{4})sG_{1}\ast G_{2}(t) \\
&=&s(1-G_{1}(t))+ (1-s) (1-\alpha )(1-h_{4})(1-G_{2}(\cdot ))\ast G_{1}(t)  +R_{i}(s)G_{1}(t)  \\
&\geq &s(1-G_{1}(t)) +sG_{1}(t)= s.
\end{eqnarray*}%
And if this is true for some $n$ then, by monotonicity
\begin{equation}
F_{n+2}^{i}(t;s)=L\left( f_{i},F_{n+1}^{i}\right) \left( t;s\right) \geq
L\left( f_{i},F_{n}^{i}\right) \left( t;s\right) =F_{n+1}^{i}(t;s)\geq s.
\label{Estss}
\end{equation}%
Next we claim that if $s\in (s_{0},1]$ then
\begin{equation}
F_{n}^{1}(t;s)\leq F_{n}^{2}(t;s), \qquad \forall n \ge0.
\label{Comp12}
\end{equation}

 Indeed, this is true for $n=0$ and if this is true for some $n$
then in view of (\ref{Estss}) for $s\in (s_{0},1]$
\begin{equation*}
F_{n+1}^{1}(t;s)=L\left( f_{1},F_{n}^{1}\right) \left( t;s\right) \leq
L\left( f_{2},F_{n}^{1}\right) \left( t;s\right) \leq L\left(
f_{2},F_{n}^{2}\right) \left( t;s\right) =F_{n+1}^{2}(t;s).
\end{equation*}%
Now on account of (\ref{Comp12}) we may pass to the limit as
$n\rightarrow \infty $ to get (\ref{CompFF}). The lemma is proved.
$\Box$

\medskip

\noindent \textbf{Proof of Theorem \ref{Tonedim4Gen}}. The first
part of the theorem is simply Proposition \ref{Lprecise}.

To prove the second part assume that $f(s)$ is not a polynomial probability
generating function (otherwise Theorem \ref{Tonedim4} gives the desired
statement). Let, for a fixed $\varepsilon >0,$ $f_{-}(s)$ and $f_{+}(s)$ be
the polynomial probability generating functions satisfying the conditions of
Lemma \ref{Lapprox} and let $F^{-}(t;s)$ and $F^{+}(t;s)$ be the probability
generating functions corresponding to the branching processes in $\mathbb{Z}%
^{4}$ with branching at the origin only and the reproduction laws
specified by $f_{-}(s)$ and $f_{+}(s)$ respectively. Let $q^{\pm
}(t):=1-F^{\pm }(t;0)$. Remark that the asymptotic of $q^{\pm}(t)$
is given by (\ref{Asympq}) with corresponding constants related to
$f^{(2)}_{\pm}(1)$.  Let $\varepsilon>0$ be small.  By
(\ref{DoublIneq}), we have that  for all sufficiently large $t$
\begin{equation*}
\frac{1}{1+2\varepsilon }\leq \frac{q^{+}(t)}{q(t)}\text{ and }\frac{q^{-}(t)%
}{q(t)}\leq \frac{1}{1-2\varepsilon }
\end{equation*}%
  while by (\ref{Asympq}), (\ref{FirstMoment}) and
  (\ref{DoublIneq}), we have that for all large $t$,
\begin{equation*}
(1-2\varepsilon) \frac{q^-(t)}{P_{1}^-(t)} \le
\frac{q(t)}{P_{1}(t)}\le  (1+2\varepsilon)
\frac{q^-(t)}{P^-_{1}(t)}, \quad (1-2\varepsilon)
\frac{q^+(t)}{P^+_{1}(t)} \le \frac{q(t)}{P_{1}(t)}\le
(1+2\varepsilon) \frac{q^+(t)}{P^+_{1}(t)},
\end{equation*}
where $P^\pm_1(t)$ are defined in the obvious way.
Clearly, for any $\lambda >0$%
\begin{eqnarray}
\mathbf{E}\left[ e^{-\lambda \mu _{0}(t)q(t)/P_{1}(t)}|\mu
_{0}(t)>0\right] &=&\frac{F\left( t;e^{-\lambda q(t)/P_{1}(t)}\right) -F\left( t;0\right) }{1-F\left( t;0\right) }
\notag \\
&=&1-\frac{1-F\left( t;e^{-\lambda q(t)/P_{1}(t)}\right) }{%
q(t)}.  \label{PreLim}
\end{eqnarray}%
Further,    $e^{-\lambda q(t)/P_{1}(t)}>s_{0}$ for all large $t$  and we deduce from   Lemma \ref{Lccomp}
that
\begin{eqnarray*}
\frac{q^{+}(t)}{q(t)}\frac{1-F^{+}\left( t;e^{-\lambda q(t)/P_{1}(t)}\right) }{q^{+}(t)}&\leq& \frac{1-F\left( t;e^{-\lambda
q(t)/P_{1}(t)}\right) }{q(t)} \\
&\leq& \frac{q^{-}(t)}{q(t)}\frac{1-F^{-}\left( t;e^{-\lambda q(t)/P_{1}(t)}\right) }{q^{-}(t)}.
\end{eqnarray*}

\noindent On the other hand, by monotonicity, $F^{+}\left(
t;e^{-\lambda q(t)/P_{1}(t)}\right)  \le
F^{+}\left( t;e^{-\lambda  (1+2\varepsilon) q^+(t)/P^+_{1}(t)}\right) $ and \linebreak
$F^{-}\left( t;e^{-\lambda
q(t)/P_{1}(t)}\right)  \ge F^{-}\left(
t;e^{-\lambda (1-2\varepsilon) q^-(t)/P^-_{1}(t)
}\right)$. Hence, letting $t\rightarrow \infty $ we get on account
of Theorem \ref{Tonedim4} that
\begin{eqnarray*}
\frac{1}{1+2\varepsilon }\left( \frac{2}{3}-\frac{2}{3}\frac{2}{2+ 3(1+2\varepsilon) \lambda }%
\right) &\leq &\lim \inf_{t\rightarrow \infty }\frac{1-F\left(
t;e^{-\lambda
q(t)/P_{1}(t)}\right) }{q(t)} \\
&\leq &\lim \sup_{t\rightarrow \infty }\frac{1-F\left( t;e^{-\lambda
q(t)/P_{1}(t)}\right) }{q(t)} \\
&\leq & \frac{1}{1-2\varepsilon }\left( \frac{2}{3}-\frac{2}{3}\frac{2}{%
2+3(1-2\varepsilon)\lambda }\right) .
\end{eqnarray*}
Letting now $\varepsilon \rightarrow +0$ we see by (\ref{PreLim}) that

\begin{equation*}
\lim_{t\rightarrow \infty }\mathbf{E}\left[ e^{-\lambda \mu
_{0}(t)q(t)/P_{1}(t)}|\mu _{0}(t)>0\right] =\frac{1}{3}+%
\frac{2}{3}\frac{2}{2+3\lambda },
\end{equation*}%
as desired. $\Box$

\begin{remark}
It follows from (\ref{FirstMoment}) and Proposition \ref{Lprecise} that the scaling in Theorem \ref{Tonedim4Gen} has the following asymptotic behavior
\begin{equation*}
\mathbf{E}\left[ \mu _{0}(t)|\mu _{0}(t)>0\right]
=\frac{\mathbf{E}\mu
_{0}(t)}{\mathbf{P}\left( \mu _{0}(t)>0\right) }\sim  {3 \over \alpha f^{(2)}(1)  C^2 } \, \frac{t}{\log ^{2}t%
},\qquad  t\rightarrow \infty .
\end{equation*}
\end{remark}

{\bf Acknowledgment } A part of this paper was written during the visit of the second author to the University Paris XIII whose support is greatly acknowledged.

\bigskip
\bigskip


{\footnotesize

\baselineskip=12pt

\noindent
\begin{tabular}{lll}
& \hskip20pt Yueyun Hu
    & \hskip60pt  Valentin Topchii \\
& \hskip20pt D\'epartement de Math\'ematiques
    & \hskip60pt Omsk branch \\
& \hskip20pt Universit\'e Paris XIII
    & \hskip60pt of Sobolev Institute of Mathematics SB RAS\\
& \hskip20pt F-93430 Villetaneuse
    & \hskip60pt Pevtcov street, 13, 644099
Omsk\\
& \hskip20pt France
    & \hskip60pt Russia \\
& \hskip20pt {\tt yueyun@math.univ-paris13.fr}
    & \hskip60pt
    {\tt topchij@ofim.oscsbras.ru}
\end{tabular}

\bigskip
\bigskip

\hskip20pt Vladimir Vatutin\par \hskip20pt Steklov Mathematical
Institute \par \hskip20pt Gubkin street 8, 119991  Moscow  \par
\par \hskip20pt Russia\par \hskip20pt {\tt vatutin@mi.ras.ru}

}

\end{document}